\documentclass[11pt]{amsart}
\usepackage{geometry}                % See geometry.pdf to learn the layout options. There are lots.
\usepackage{graphicx}
\usepackage{amssymb}
\usepackage{epstopdf}
\usepackage{amsmath}
\usepackage{amsfonts}
\usepackage{amssymb}
\usepackage{amsthm}

\usepackage[parfill]{parskip}

\usepackage{mathrsfs}
\usepackage{latexsym}
\usepackage{longtable}    % may not have to use this package
\usepackage{graphpap}     % may not have to use this package
\DeclareMathAlphabet{\mathpzc}{OT1}{pzc}{m}{it}

\newtheorem{thm}{Theorem}[section]
\newtheorem{lemma}[thm]{Lemma}
\newtheorem{cor}[thm]{Corollary}
\newtheorem{prop}[thm]{Proposition}
\theoremstyle{definition}
\newtheorem{defin}[thm]{Definition}
\newtheorem{ex}[thm]{Example}
\theoremstyle{remark}
\newtheorem{rem}[thm]{Remark}
\numberwithin{equation}{section}
\numberwithin{thm}{section}

\def \e {\varepsilon}
\def \F {\mathcal F}

\def \R {\mathbb R}
\def \Z {\mathbb Z}
\def\Q {\mathbb Q} 
\def\N {\mathbb N}
\def\C {\mathbb C}

\newcommand{\kf}{\mathbb{K}}

\newcommand{\lm}{\lambda}
\newcommand{\lid}{\lambda I - D}
\newcommand{\eo}{E_{\omega}}
\newcommand{\se}{\sigma_{e}}
\newcommand{\spp}{\sigma_{p}}
\newcommand{\sep}{\sigma_{e}^{\prime}}
\newcommand{\sepp}{\sigma_{e}^{\prime\prime}}
\newcommand{\plam}{\partial \Lambda}
\newcommand{\knav}{(\kf, |\cdot|)}
\newcommand{\kft}{\kf^{t}}
\newcommand{\vp}{V_{p}}
\newcommand{\qp}{\Q_{p}}
\begin{document}

\title{Spectral Analysis for a class of Bounded Linear Operators}
\author{Teylama Herve Miabey}
\address{Department of Mathematics \\ University of the District of Columbia} 
\date{}                                           % Activate to display a given date or no date
\maketitle

\begin{abstract}
\thispagestyle{plain}
\setcounter{page}{6}       %this number may need to be changed
\addcontentsline{toc}{chapter}{Abstract}
\indent
We study the Spectral Analysis for
a class of bounded linear operators $T = D + u_{1} \otimes v_{1} +
u_{2} \otimes v_{2} + \cdots + u_{m} \otimes v_{m}$ in a non 
Archimedean Hilbert space $\eo$, where $D$ is a diagonal linear 
operator and $F = u_{1} \otimes v_{1} +
u_{2} \otimes v_{2} + \cdots + u_{m} \otimes v_{m}$
is a finite rank linear operator. In fact, our study is a 
generalization of the work of Diagana and McNeal. In their joint paper,
Diagana and McNeal study the Spectral Analysis for rank 1 perturbation
of a diagonal operator in $\eo$ defined by $A = D_{\lambda} +
u \otimes v$. Under some suitable assumption, they prove that 
\[
\sigma(A) = \{\theta_{j} \}_{j\ge 0}\cup \spp(A)
\]
where $\spp(A)$ is the set of all the eigenvalues of $A$ and 
$\theta = \left( \theta_{j}\right)_{j\in \N}$, $\theta_{j} =
\lambda_{j} + \omega_{j}\alpha_{j}\beta_{j}$, with all $j \in \N$. 

In this study of the Spectral Analysis, we use extensively the Theory
of Fredholm Operators to deduce some of our main results. In a forthcoming paper,
we will apply these methods to study a similar class of unbounded linear operators.
\end{abstract}

\section{Introduction and Background}
\begin{defin}
Let $\kf$ be an arbitrary field. We define an \emph{absolute value}
on $\kf$ to be the function
\[
|\cdot| : \kf \rightarrow  [0,+\infty)
\]
subject to the following conditions:
\begin{enumerate}
\item[(1)] $|x|= 0$ iff $x=0 \quad \forall x \in \kf$
\item[(2)] $|xy| = |x||y| \quad \forall x,y \in \kf$
\item[(3)] $|x+y| \le |x| +|y|\quad \forall x,y \in \kf$. 
\end{enumerate}
\end{defin}

$\kf$ endowed with the aforementioned absolute value 
$(\kf, |\cdot|)$ is called a \emph{non Archimedean value field}
 if condition (3) is replaced by the non Archimedean triangle
inequality 
\begin{enumerate}
\item[($3^{\prime}$)]$|x+y| \le \max\{|x|,|y|\}$.  
\end{enumerate}

If $|x| \ne |y|$, then $|x+y| = \max \{|x|,|y|\}$.

\begin{rem}
Note that $|x+y| \le \max \{|x|,|y|\} \le |x| + |y|$, which implies
that condition $(3^{\prime})$ is stronger than condition $(3)$. 
\end{rem}
 
\begin{defin}
Let $\kf$ be a field. A non Archimedean metric $d$ on $\kf$ is a 
function
\[
d : \kf \times \kf \rightarrow \R
\]

such that 
\begin{enumerate}
\item[(1)] $d(x,y) \ge 0$ for all $x,y \in \kf$
\item[(2)] $d(x,y) = 0$ iff $x = y$
\item[(3)] $d(x,y) = d(y,x)$ for all $x,y \in \kf$
\item[(4)] $d(x,y) \le max\left\{ d(x,z), d(z,y)\right\}$ 
for all $x,y$ and $z \in \kf$. 
\end{enumerate}

The pair $(\kf, d)$ is called a non Archimedean metric space. 

\end{defin}

\begin{defin}
Let $\kf$ be a field and $|\cdot|$ be a non Archimedean absolute
value on $\kf$. The non Archimedean distance between two elements
$x,y$ of $\kf$ is defined by $d(x,y) = |x-y|$. 
\end{defin}

\begin{rem}
\begin{enumerate}
\item[(1)] If $|\cdot|$ is a non Archimedean absolute value on $\kf$, then 
$d(x,y) = |x-y|$ is a non Archimedean metric on $\kf$. 
\item[(2)] Clearly, the non Archemedean absolute value induces a 
meteoric $d$ on $\kf$ defined above. Therefore, we can consider 
the notion of convergence in $\kf$. A sequence 
$(x_{n})_{n \in \N} \subset \kf$ is said to converge to $x \in \kf$ if
$d(x_{n}, x) \rightarrow 0$ as $n \rightarrow \infty$. The 
non Archimedean value field $(\kf, |\cdot|)$ is complete if the
metric space $(\kf, d)$ is complete. $(\kf, d)$ is complete iff
all the Cauchy sequences in $\kf$ converge in $\kf$. 
\end{enumerate}

\end{rem}

\begin{prop}
Let $\kf$ be a field and $|\cdot |$ be a normed value field. Then the 
following conditions are equivalent:
\begin{enumerate}
\item[(1)] $|\cdot|$ is a non Archimedean absolute value
\item[(2)] $|\cdot|$ is bounded in $\Z$. 
\item[(3)] $|x|\le 1$ for all $x \in \Z$. 
\end{enumerate}
\end{prop}

\begin{proof}
$(1) \Rightarrow (3)$:

Note that since $|\cdot|$ is non Archimedean, it implies that 
$|x| \le 1$ since for any $x \in \Z$

\[
\begin{aligned}
|x| & = |x\cdot 1| = |\underbrace{1+1+\cdots +1}_{x \text{ times}}|\\
&\le \max\{|1|, |1|,\dots,|1|\} \\
& \Rightarrow |x| \le 1, \text{ since } \max\{|1|, |1|,\dots,|1|\} = 1
\end{aligned}
\]

$(3) \Rightarrow (2)$:

Since $|x| \le 1$ for any $x\in \Z$, it implies that $|\cdot|$ is 
bounded in $\Z$. Therefore to complete the proof, it is sufficient
to show that if $|\cdot|$ is bounded in $\Z$, then $|\cdot|$ is a 
non Archimedean absolute value, which is $(2) \Rightarrow (1)$. 

Suppose that $|m| \le M$ for all $m \in Z$. If $x,y \in \kf$ with 
$n \in \Z$, then 
\[
\begin{aligned}
|x+y|^{n} &= |(x+y)^{n}| = \left|\sum_{i=0}^{n}\binom{n}{i} x^{i}y^{n-i}  
\right|\\
&\le \sum_{i=0}^{n}\left|\binom{n}{i}\right||x|^{i}|y|^{n-i}
\quad \quad \left(\text{Note that $\binom{n}{i} \le M$}.\right)\\
&\le M \sum_{i=0}^{n}(\max(|x|,|y|)^{n}=(n+1)M\max\left\{|x|,|y|\right\}^{n}\\
\text{Therefore, } |x+y|^{n} &\le (n+1)M\max\left\{|x|,|y|\right\}^{n}\\
\text{Hence, } \sqrt[n]{|(x+y)^{n}}&\le (n+1)^{1/n}M^{1/n}\max\left\{
|x|,|y|\right\} \\
|x+y| &\le \sqrt[n]{(n+1)M} \max\left\{|x|,|y|\right\} \\
\text{but } &\lim_{n\rightarrow \infty}\sqrt[n]{(n+1)M}= 1\\
\text{Therefore, }|x+y| &\le \max\left\{|x|,|y|\right\}. \\
\end{aligned}
\]
\end{proof}

\begin{rem}
The difference between the Archimedean and non Archimedian Absolute
values can be established in the following way: An Archimedean absolute
value implies that  for every $x,y \in \kf$ such that $x\ne 0$, 
there exists $n \in \Z^{+}$ such that $|nx| > |y|$, i.e. 
$\sup \left\{\,|n|: n \in \Z \,\right\} = +\infty$. However the non 
Archimedean absolute value implies that 
$\sup\{\,|n| : n \in \Z\,\}=1$.
\end{rem}

\begin{prop}
Let $(\kf, |\cdot|)$ be a non Archimedean value field. Let 
$|x| \ne |y|$, then $|x+y = \max\left\{|x|,|y|\right\}.$
\end{prop}

\begin{proof}
Let $x,y \in \kf$ and suppose that $|x| > |y|$. Then
\[
\begin{aligned}
|x+y|\le \max \left\{|x|,|y|\right\} \\
& \le |x| = |(x+y) - y| \le \max\left\{|x+y|,|y|\right\} \\
\text{And since} |y| < |x|\\
\max\left\{|x+y|,|y|\right\} &= |x+y|\\
\text{but } |x+y| &\le \max\left\{x|,|y|\right\} = |x|\\
\end{aligned}
\]
We have thus established that $|x+y|\le |x|$ and $|x| \le |x+y|$. Then
it is clear that $|x+y| = |x|$, $\Rightarrow$ if $|x| \ne |y|$
then $|x+y| = \max\left\{|x|,|y|\right\}.$
\end{proof}

\begin{prop}
Let $(\kf, |\cdot|)$ be a metric space on a non Archimedean value
field and let $x, y $ and $z\in \kf$ such that $|x-z| \ne |z-y|$. Then
$|x-y|= \max\left\{|x-z|, |z-y|\right\}$
\end{prop}

\begin{proof}
Suppose that $|z-y| < |x-z|$, then we have 
$|x-y| \le \max\left\{|x-z|, |z-y|\right\} = |x-z|$ and suppose now
that $|x-y| <|x-z|$. Consequently,
\[
\begin{aligned}
|x-z| &\le \max \left\{|x-y|, |y-z| \right\} < |x-z|\\
\text {which gives}\\
|x-z| &< |x-z|, \quad \text{ which is a contradiction}\\
\therefore \quad |x-z| &= \max\left\{|x-z|, |z-y|\right\}\\
\end{aligned}
\]
\end{proof}

\begin{cor}
In the non Archimedean space, all triangles are isosceles.
\end{cor}

\begin{proof}
Consider $x,y$ and $z$ to be the vertices of a given triangle. Note 
that $|x-y| = |(x-z)+(z-y)|$. Suppose now $|x-z| = |z-y|$. Then 
two of the sides of the triangle have equal lengths, which implies that
the triangle is isosceles. Suppose now $|x-z| \ne |z-y|$ and that 
$|x-z| > |z-y|$. 

Since we are in a non Archimedian field, we have $|x-y \le
\max\left\{|x-z|,|z-y|\right\}$ with equality when $|x-z| \ne |z-y|$. 
Therefore $|x-y|  = \max\left\{|x-z|,|z-y|\right\} = |x-z|$. Therefore
$|x-y| = |x-z|$ and again two of the sides of the triangle have
equal lengths. Therefore all triangles in the  non Archimedean 
field are isosceles.
\end{proof}

\begin{defin}
Let $(\kf, |\cdot|)$ be a non Archimedean value field. 
\begin{enumerate}
\item[(1)] A sequence $(x_{n})_{n\in \N} \subset \kf$ converges to
$x \in \kf$ if for each $\e >0$, $\exists\, \eta_{0} \in \N$ such that
$|x_{n}-x| < \e$ whenever $n\ge \eta_{0}$. 
\item[(2)] A sequence $(x_{n})_{n\in \N} \subset \kf$ is called a 
Cauchy sequence if for every $\e > 0$, $\exists\, \eta_{0} \in \N$ 
such that $|x_{n} - x_{m}| < \e$ whenever $m\ge n\ge \eta_{0}$. 
\end{enumerate}
\end{defin}

Due to the non Archimedian triangle inequality, we conclude the
following proposition:

\begin{prop}
If $(\kf, | \cdot |)$ is a non Archimedean value field, 

a sequence
$(x_{n})_{n\in \N} \subset \kf$  is a Cauchy sequence iff 
\[
|x_{n+1} - x_{n}| \rightarrow 0 \text{ as } n \rightarrow \infty.
\] 
\end{prop}

\begin{proof}
$(\Rightarrow)$:

Suppose that $(x_{n})_{n\in \N}$ is a Cauchy sequence, then by the 
definition of the Cauchy sequence, for every $\e > 0$, there exists
a $\eta_{0} \in \N$ such that $|x_{n} - x_{m}| < \e$ whenever 
$m \ge n \ge \eta_{0}$ which implies that $|x_{n}-x_{m}| \rightarrow 0$
as $n \rightarrow \infty$.

$(\Leftarrow)$:
Suppose that $|x_{n+1}-x_{n}|\rightarrow 0$ as $n \rightarrow \infty$
and suppose that $m\ge n$. We have 

$x_{m}-x_{n}=x_{m}-x_{m-1}+x_{m-1}-x_{m-2}+x_{m-2}-\cdots+x_{n+1}-x_{n}$ .
Combining each pair of terms of the form $x_{t}-x_{t-1}$, 
$t = m,\dots, n$ and based on the earlier definition, we conclude 
that $|x_{m} -x_{n}| \rightarrow 0$ as $n, m \rightarrow \infty$. 
Therefore $(x_{n})_{n \in \N}$ is a Cauchy sequence in $\kf$. 
\end{proof}

\begin{prop}
Let $(\kf, |\cdot|, d)$ be a complete non Archimedean field. The series
$\displaystyle \sum_{n=0}^{\infty}x_{n}$ converges to $x\in \kf$
iff $x_{n}\rightarrow 0$ as $n \rightarrow \infty$.
\end{prop}

\begin{proof}
Set $\displaystyle S_{n} = \sum_{r=0}^{n}x_{r}$. It follows that
$x_{n}= S_{n} - S_{n-1}$. Then if $S_{n} \rightarrow S$, then
$x_{n}\rightarrow 0 $ in $\kf$ as $n \rightarrow \infty$. Conversely,
if $x_{n}\rightarrow 0$, then $(S_{n})_{n\in \N}$ is a Cauchy sequence
by the previous proposition. We know that $\kf$ is complete, hence 
there exists an element $s \in \kf$ such that $S_{n}\rightarrow S$
as $n\rightarrow \infty$, which completes the proof.
\end{proof}

\begin{defin}
Let $(\kf, |\cdot|)$ be a non Archimedean value field. 
\begin{enumerate}
\item[(1)] An open ball of radius $r$ and center $x$ is defined by

$B(x,r) = \{y \in \kf:|x-y|<r\}$.

\item[(2)] A closed ball of radius $r$ and center $x$ is defined by

$\overline{B(x,r)} = \{y \in \kf:|x-y|\le r\}$.

\item[(3)] A sphere in $\kf$ of radius $r$ and center $x$ is defined
by $S(x,r) = \{y \in \kf:|x-y| = r\}$.
\end{enumerate}
\end{defin}

\begin{prop}
Let $\knav$ be a non Archimedean value field. Then
\begin{enumerate}
\item[(1)] Every point that is contained in an open or closed ball
in $\knav$ ia a center of the ball.
\item[(2)] Every ball is both open and closed. These balls are called
clopen.
\item[(3)] If $a,b$ are elements of $\kf$ and $r_{1}, r_{2} \in
\R_{+}^{\times}$, then $B(a,r_{1}) \cap B(b,r_{2}) \ne \phi$ iff 
$B(a,r_{1}) \subset B(b,r_{2})$ or $B(a,r_{1}) \supset B(b,r_{2})$, 
which implies that any two balls are either disjoint or one ball 
contains another. 
\end{enumerate}
\end{prop}

\begin{proof}
For a proof of (3), see Gouvea \cite{Gou1} pg 34.
\begin{enumerate}
\item[(1)] Let $B(a,r)$ be a ball in $\kf$ and let $y \in B(a,r)$ and
$z \in B(y,r)$. Therefore
\[
\begin{aligned}
|z-a| &= |z-y+y-a|\le \max\{|z-y|,|y-a|\}\\
|z-a| &\le \max\{|z-y|,|y-a|\} < r\\
\end{aligned}
\]
Therefore $z \in B(a,r) \Rightarrow  B(y,r) \subset B(a,r)$.
\[
\begin{aligned}
w \in B(a,r) &\Rightarrow |w-y| = |w-a+a-y| \\
& \le \max\{|w-a|,|a-y|\} <r \\
\Rightarrow B(a,r) \subset B(y,r) &\Rightarrow B(y,r) = B(a,r)\\
\end{aligned}
\]
Therefore every point in the closed or open ball is a center of the
ball.
\item[(2)] Let $B(a,r)$ be a ball in $\knav$. $z \in B(a,r),
\exists \varepsilon > 0, \exists w \in B(a,r) \cap B(z,r)$.

Let $\e \le r \Rightarrow |w-a| < r, |w-z| < \e$

$\Rightarrow |z-a| = |z-w+w-a|\le \max\{|z-w|, |w-a|\}<r$.

Since $|z-w| < \e$ and $|w-a| < r$, $\Rightarrow z \in B(a,r)$.
Let $\overline{B(a,r)}$ be a closed ball, $w \in \overline{B(a,r)}$ and
let $\e \le r$. Then 
\[
\begin{aligned}
z \in B(w,\e)  \Rightarrow |z-a| &= |z-w+w-a| \\
|z-a| \le \max\{|z-w|, |w-a|\} <r.
\end{aligned}
\] 
Since $|z-w|<\e$ and $|w-a| < r$, then $|z-a|\le r \Rightarrow
z \in \overline{B(a,r)}$.
\end{enumerate}
\end{proof}

\begin{prop}
Let $\knav$ be a non Archimedean field. Then
\begin{enumerate}
\item[(1)] $\overline{B(0,1)} = \{ x \in \kf: |x|\le 1\}$
is a subring of $\kf$. 
\item[(2)] $B(0,1) =\{ x \in \kf: |x| < 1\}$
is an ideal of $\overline{B(0,1)}$. 
\item[(3)] $B(0,1)$ is a unique maximal ideal of $\overline{B(0,1)}$.
\item[(4)] The Quotient $\displaystyle \overline{B(0,1)}/B(0,1)$
is a field. It is the residue class field of $|\cdot|$. 
\item[(5)] $\kf$ is a field of fractions of $\overline{B(0,1)}$. 
\end{enumerate}
\end{prop}

\begin{proof}
See Gouvea \cite{Gou1}.
\end{proof}

\section{The Vector Space $\kf^{t}$}
\begin{defin}
Let $\knav$ be a non Archimedean field. Just as in the classical 
case of $\R$ and $\C$, there is a vector space $\kft$ which is 
associated to $\kf$ and defined as follows:
\[
\kft = \left\{\, x= (x_{1},x_{2}, \dots, x_{t}): x_{r}\in \kf, 
r = 1, \dots, t\,\right\}.
\]
Note that the norm defined by 
$\displaystyle ||x||_{t} = \max_{1\le r\le t}|x_{r}|$ 
for all $x = (x_{1},x_{2}, \dots, x_{t}) \in \kft$ is a non
Archimedean norm over $\kft$. The inner product on $\kft$ is defined
by 
\[
\langle x,y\rangle_{t} = \sum_{r=1}^{t} x_{r}y_{r}=
x_{1}y_{1}+x_{2}y_{2}+ \cdots + x_{t}y_{t}
\]
where $x = (x_{1},x_{2}, \dots, x_{t})$ and $y = (y_{1},y_{2}, \dots, y_{t}) \in \kft$. 

The inner product satisfies the Cauchy-Schwartz inequality i.e.
$|\langle x,y\rangle_{t}| \le ||x||_{t}||y||_{t}, x,y \in \kft$, and 
any element $x = (x_{1},x_{2}, \dots, x_{t})$ can be written 
uniquely as $\displaystyle x = \sum_{r=1}^{t} x_{r}e_{r}$, where
$e_{r}= (0,0,\dots,1,0,\dots,0)$ and $(e_{r})_{r=1,\dots,t})$
is the orthonormal base of $\kft$. 
\end{defin}

\section{A Valuation}
\begin{defin}\label{def-valuation}
A valuation on an arbitrary field $\kf$ is a function 
$V :\kf \rightarrow \R \cup \{\infty\}$ such that for every 
$x,y\in \kf$, we have the following properties:
\begin{enumerate}
\item[(1)] $V(x) = \infty$ iff $x = 0$
\item[(2)] $V(xy) = V(x) + V(y)$
\item[(3)] $V(x+y) \ge \min \{\,V(x),V(y)\,\}$
\end{enumerate}
\end{defin}

\begin{rem}
For any element $c \in \kf$ such that $c > 1$, $c^{-V(\cdot)} = 
|x|$ defines a non Archimedean absolute value on $\kf$. The valuation
is called non trivial if there exists an element $a \in \kf^{*}$
such that $V(a) \ne 0$. 
\end{rem}

\begin{defin}
A valuation $V$ on $\kf$ is said to be a discrete valuation if
the totally ordered group $V(\kf^{*})$ is isomorphic to $\Z$. 
\end{defin}

\begin{lemma}
Let $V$ be a discrete valuation, $x,y \in \kf$, $V(x) \ne V(y)$. 
Then we have 
\[
V(x+y) = \min\{\,V(x), V(y)\,\}
\]
\end{lemma}

\begin{proof}
Since $V(x) \ne V(y)$, suppose without loss of generality that 
$V(x) < V(y)$. Suppose now that $V(x+y) \ne \min \{\,V(x), V(y)\,\}$.
We have $V(x+y) > V(x)$. 

Note that $V(x) = V[(x+y)-y] \ge
\min\{V(x+y),V(y)\}$, since $V(-y) = V(y)$. Therefore,
$V(x) = V[(x+y)+y]\ge \min\{V(x+y),V(y)\} > V(x)$. Hence 

$V(x)
> V(x)$, a contradiction. Therefore $V(x+y) = \min\{V(x),V(y)\}$.  
\end{proof}

\begin{ex}
Let $\F = \kf(x)$. Define $V_{p(x)}: \F \rightarrow \Z \cup \{\infty
\}$, satisfying the conditions of Definition \ref{def-valuation}, 
where $p(x) \in \kf[x]$ is an irreducible monic polynomial. $V_{p(x)}$
is a discrete valuation. We define the discrete valuation ring as the
set \[
\mathcal{O}_{p(x)} = \left\{ \frac{f(x)}{g(x)}: f(x), g(x)
\in \kf[x](p(x),g(x)) = 1 \right\}
\]
and its unique maximal ideal 
\[
P_{p(x)} = \left\{ \frac{f(x)}{g(x)}: f(x), g(x)
\in \kf[x], p(x)|f(x) \text{ and } p(x)\nmid g(x) \right\}.
\]
There is another valuation ring on $\kf(x)/\kf$:
\[
\mathcal{O}_{\infty} = \left\{ \frac{f(x)}{g(x)}: f(x), g(x)
\in \kf[x], \deg f(x) \le \deg g(x) \right\}
\]
and its maximal ideal 
\[
P_{\infty} = \left\{ \frac{f(x)}{g(x)}: f(x), g(x)
\in \kf[x], \deg f(x) < \deg g(x) \right\}.
\]
\end{ex}

\section{The Field $\mathbb{Q}_{p}$}

\begin{defin}
A p-adic valuation on $\Z$ is a function $V_{p}: \Z-\{0\}
\rightarrow \R$ where $p \in \Z$ is a prime number and $V_{p}(t)$ is
a unique positive number such that $t = p^{V_{p}(t)}\cdot q$ and
$(p,q) = 1$. 
\end{defin}

\begin{rem}
The valuation $V_{p}$ on $\Z$ can be extended to the field of rational
numbers $\Q$, i.e. if $x \in \Q - \{0\}$, then $x = t/s$, where
$(t,s) \in \Z^{2}-\{0\}$. 

$\displaystyle V_{p}(x) = V_{p}\left(\frac{t}{s}\right)=
V_{p}\left(t\cdot \frac{1}{s}\right) = V_{p}(t) - V_{p}(s)$. If $x = 0$, 
we have $V_{p}(0) = \infty$. 
\end{rem}

\begin{cor}
Consider $(x,y) \in \Z^{2}$. We have the following properties
of the valuation:
\begin{enumerate}
\item[(1)] $\vp(x) = \infty$ iff $x = 0$.
\item[(2)] $\vp(x\cdot y) = \vp(x) + \vp(y)$
\item[(3)] $\vp(x+y) \ge \min \left\{\vp(x), \vp(y) \right\}$, with 
equality if $\vp(x) \ne \vp(y)$. 
\end{enumerate}
\end{cor}

\begin{proof}
\begin{enumerate}
\item[(1)] is deduced from the definition.
\item[(2)] Let $x,y \in \Z$. Then $x = p^{r}x^{\prime}$, 
$y = p^{s}y^{\prime}$, where $x^{\prime}, y^{\prime}$ are integers
and $p$ is relatively prime with respect to $x^{\prime}$ and 
$y^{\prime}$ and $r, s \in \Z$. Suppose that $r > s$. 
\[
\begin{aligned}
x \cdot y & = p^{r}\cdot x^{\prime}\cdot p^{s}\cdot y^{\prime} \\
& = p^{r}\cdot p^{s}\cdot x^{\prime}y^{\prime}\\
&= p^{r+s}\cdot x^{\prime}y^{\prime}.\\
\vp(x,y) &= \vp(p^{r+s}\cdot x^{\prime}y^{\prime})\\
&= \vp(p^{r+s}) + \vp(x^{\prime}y^{\prime})\\
&= r+s, \quad \text{ since } \vp(x) = r \text{ and } \vp(y) = s\\
\therefore \vp(x \cdot y) &= \vp(x) + \vp(y) \\
\end{aligned}
\]
\item[(3)] In this case, suppose that $x$ and $y$ are rational
numbers. Then we have $\displaystyle x = p^{r}\frac{a}{b}$ and
$\displaystyle y = p^{s}\frac{c}{d}$, where $a, b, c$ and $d \in \Z$,
$p$ is relatively prime with respect to $a, b, c$ and $d$ and 
$r,s \in \Z$. Now suppose that $r = s$. Then 
\[
\begin{aligned}
x+y &= p^{r}\left(\frac{a}{b} + \frac{c}{d} \right) \\
&= p^{r}\cdot \left(\frac{ad+bc}{bd} \right) \\
\Rightarrow \vp(x+y) &= \vp\left[p^{r}\left( \frac{ad+bc}{bd}
\right) \right] \\
\Rightarrow \vp(x+y) &\ge r \quad \text{ since } (p,bd) = 1. \\
\end{aligned}
\]  
Assume that $r \ne s$ and set $s > r$. Then we have 
\[
\begin{aligned}
x+y &= p^{r}\frac{a}{b} + p^{s}\frac{c}{d} \\
&= p^{r}\left(\frac{a}{b} + p^{s-r}\cdot \frac{c}{d} \right) \\
&= p^{r}\left( \frac{ad+bc\cdot p^{s-r}}{bd}\right). \\
\end{aligned}
\]
Note that since $s>r \Rightarrow s-r>0$ and since $(p,bd)=1
\Rightarrow \vp(x+y) = r = \min\{\vp(x), \vp(y) \}$. 
\end{enumerate}
\end{proof}

\begin{defin}
For $x \in \Q$, the p-adic absolute value of $x$ is given by 
\[
|x|_{p}= \left\{\begin{array}{cc}
p^{-\vp(x)} & \text{ if } x \ne 0 \\
p^{-\infty}= 0 & \text{ if } x = 0 \\
\end{array}
\right.
\]
where $\vp(s)$ is a valuation of $x\in \Q$. 
\end{defin}

\begin{rem}
The p-adic absolute value is a non Archimedean valuation.
\end{rem}

\begin{prop}\label{navfield}
Consider the function $|\cdot |_{p}: \Q \rightarrow \R^{+} \cup \{0\}$.
Then we have the following properties:
\begin{enumerate}
\item[(1)] $|x|_{p} = 0 $ iff $x = 0$. 
\item[(2)] $|xy|_{p} = |x|_{p}|y|_{p}$. 
\item[(3)] $|x+y|_{p} \le \max\{|x|_{p}, |y|_{p}\}$. 
\item[(4)] $|x|_{p}\ne |y|_{p}\Rightarrow |x+y|_{p} = 
\max\{|x|_{p}, |y|_{p}\}$. 
\end{enumerate}
Hence $(\Q, |\cdot|_{p})$ is a non Archimedean value field. 
\end{prop}

\begin{rem}
The p-adic absolute value $|\cdot|_{p}$ induces a metric 
$d$ on $\Q \times \Q$ defined by $d(x,y) = |x-y|$ for all $x,y \in 
\Q$. 
\end{rem}

\textsc{Proof of Proposition \eqref{navfield}}
\begin{enumerate}
\item[(1)] $|x|_{p} = 0 = p^{-\vp(x)}\Rightarrow \vp(x) = 
+\infty \Rightarrow x = 0$.
\item[(2)]\[
\begin{aligned}
|xy|_{p} &= |p^{\vp(x)}\cdot q \cdot p^{\vp(y)}\cdot q^{\prime}|_{p}\\
&= |p^{\vp(x)+\vp(y)}|_{p} = p^{-(\vp(x)+\vp(y))} \\
&= p^{-\vp(x)}\cdot p^{-\vp(y)} \\
&= |x|_{p}\cdot|y|_{p} \\
\end{aligned}
\]
\item[(3)] \[
\begin{aligned}
|x+y|_{p} &= |p^{\vp(x)}\cdot q + p^{\vp(y)}\cdot q^{\prime}|_{p}\\
&\text{Note that by definition, }\\
|x+y|_{p} &= p^{-\vp(x+y)} \\
\text{but } p^{\vp(x+y)} &\le p^{-\min\{\vp(x),\vp(y) \}}\\
&= \max\left\{p^{-\vp(x)}, p^{-\vp(y)}  \right\}\\
&= \max \left\{|x|_{p}, |y|_{p}  \right\}\\
\end{aligned}
\]
Hence $|x+y| \le \max\left\{|x|_{p}, |y|_{p}  \right\}$.
\end{enumerate}

\section{Construction of $\Q_{p}$}
\begin{prop}
For all primes $p \in \Z$, the field $\Q$ is not complete with respect
to the non Archimedean absolute value $|\cdot |_{p}$. 
\end{prop}

\begin{proof}
See Gouvea\cite{Gou1}, page 50.
\end{proof}

\begin{defin}
The field of p-adic numbers $\Q_{p}$ is the completion of the
non Archimedean value field $(\Q, |\cdot|_{p})$. 
\end{defin}

\begin{lemma}
Let $|\cdot|_{p}$ be a non Archimedean absolute value on $\Q$. Then
the extension of $|\cdot|_{p}$ to $\qp$ is also a non Archimedean
absolute value.
\end{lemma}
\begin{proof}
Let $a,b \in \qp$ and let $\displaystyle\{a_{t}\}_{t\in \N}$ and 
$\displaystyle\{b_{t}\}_{t\in \N}$ be sequences that converge to 
$a$ and $b$ respectively, i.e. $\displaystyle \lim_{t\rightarrow 
\infty}\{a_{t}\} = a$ and 
$\displaystyle \lim_{t\rightarrow 
\infty}\{b_{t}\} = b$. We then have 
\[
\begin{aligned}
 \lim_{t\rightarrow \infty} |a_{t} + b_{t}|_{p}
& \le \max\left(|a|_{t}, |b|_{t}\right) \\
& \Rightarrow |a+b| \le \max \left\{|a|,|b|\right\} \\
\end{aligned} 
\]
which concludes the proof. 
\end{proof} 

\begin{rem}
There are several methods that are used to complete 
$(\Q, |\cdot|_{p})$. We will cover three methods of completion:
the analytic method, the semi algebraic method and the algebraic 
method. 
\begin{enumerate}
\item[(1)] Analytic Method:

Let $\mathcal{C}$ be the set of all Cauchy sequences in $\Q$ with 
respect to the Absolute value $|\cdot|_{p}$. We define the following 
relation: Two sequences $\{x_{n}\}$ and $\{y_{n}\}$ of $\mathcal{C}$
are related, i.e. $\{x_{n}\} \sim \{y_{n}\}$ iff $\forall \varepsilon
> 0, \exists N > 0$ such that $n \ge N \Rightarrow 
|x_{n} - y_{n}| < \varepsilon$. 

We can easily verify that $\sim$ is an equivalence relation. We 
define $\qp = \mathcal{C}/\sim$, where addition is defined by 
$\displaystyle \overline{\{x_{n}\}} + \overline{\{y_{n}\}}
= \overline{\{x_{n}\}+\{y_{n}\}}$ and multiplication by 
$\displaystyle \overline{\{x_{n}\}} \cdot \overline{\{y_{n}\}}
= \overline{\{x_{n}\}\cdot\{y_{n}\}}$.

Therefore the elements of $\qp$ are the equivalence classes of the
Cauchy sequences. 
\end{enumerate}
\end{rem}

\section{Computation of the Essential Spectrum of $D$, $\se{D}$ where
$D$ is a bounded linear operator}
\begin{defin}
Let $\kf$ be a non Archimedian field and let $E$ and $F$ be vector spaces over
$\kf$ such that $f : E \longrightarrow F$ is a linear map. Let $\lambda \in
\kf$. We define the \emph{weak multiplicity} $\lambda$ relative to $f$ to 
be the dimension of the kernel $N(\lambda I - D)$. The multiplicity
$\lambda$ is said to be finite if the dimension of the kernel $N(\lambda I - D)$
is finite. 
The multiplicity $\lambda$ is said to be infinite if the dimension of
$N(\lambda I - D)$ is infinite. We denote $\eta(\lambda I - D)$ to be the
dimension of the kernel of $\lambda I  - D$. 
\end{defin}

\begin{thm}
Let $D: E_{\omega} \longrightarrow E_{\omega}$ be a bounded diagonal operator
such that $D(u) = \sum \lambda_{i}u_{i}e_{i}$ with $u = \sum u_{i}e_{i}$
for $\lambda_{i} \in \kf$, $\forall i \in \mathbb{N}$. Let $\lambda$ be an
eigenvalue of $D$ with finite multiplicity with respect to $D$. Then 
$\lambda I - D$ is a Fredholm operator of index 0.
\end{thm}
\begin{proof}

1. We want to show that $\eta(\lambda I-D)$ is equal to the number of 
$\lambda_{i} \in \Lambda = \left\{\,\lambda_{k} \,\right\}_{k\in \mathbb{N}}$
such that $\lambda_{k} = \lambda$.

\begin{displaymath}
\begin{aligned}
\text{Let }u \in N(\lambda I-D) & \Rightarrow (\lambda I-D)(u) = 0\\
&\Rightarrow \sum(\lambda u_{i}- \lambda_{i}u_{i})e_{i}= 0\\
&\Rightarrow \sum(\lambda - \lambda_{i})u_{i}e_{i} = 0\\
&\Leftrightarrow (\lm - \lm_{i})u_{i}= 0\quad \forall i \in \mathbb{N}\\
\end{aligned}
\end{displaymath}
For $k \in \kf$, if $\lm_{k}\ne \lm $, $u_{k} = 0$. This implies that
$N(\lm  I -D)$ is generated by $e_{i}$ such that 
$\lm_{i} = \lm \Rightarrow \eta(\lm I -D) = \text{Cardinal}\{e_{i}|\lm_{i}=\lm\}$

2. We want to compute the dimension of the cokernel  of $(\lm I-D)$, denoted
by $\delta(\lm I-D)$. Let $v \in Im(\lm I-D) \Rightarrow \exists u \in 
E_{\omega}$ such that
\begin{displaymath}
\begin{aligned}
(\lm I - D)(u) = v &
\Leftrightarrow \sum_{i\in \mathbb{N}}(\lm - \lm_{i})u_{i}e_{i}
=\sum_{i\in \mathbb{N}}v_{i}e_{i}\\
&\Leftrightarrow (\lm - \lm_{i})u_{i}=v_{i}\, \text{ if } \lm = \lm_{i}
\Rightarrow v_{i}= 0 \, \forall i \in \mathbb{N} \\
\end{aligned}
\end{displaymath}
This implies $Im(\lambda I  - D)$ is generated by the $e_{k}$ such
that $\lambda_{k} \ne \lambda \Rightarrow E_{\omega}/Im(\lambda I - D)$
is of the dimension equal to the number of $\{\, e_{i}| 
\lambda_{i}= \lambda \,\}$ or $dim(E_{\omega}/Im(\lambda I - D))
= \delta(\lambda I - D)$. Therefore $\delta(\lambda I  - D)
= \# \{\, e_{i}| \lambda_{i}= \lambda \,\}$

3. We want to show that if $\lambda$ is an eigenvalue with finite
multiplicity with respect to $D$, then the range, denoted by 
$R(\lid)$ is closed in $E_{\omega}$. Let $F$ be a vector subspace of
$E_{\omega}$. We define the orthogonal set of $F$, denoted 
$F^{\perp}$ to be the set of all $v \in E_{\omega}$ such that
$(u,v) = 0\, \forall u \in F$, i.e. $F^{\perp} = \{\, v \in E_{\omega}
\,|\, (u,v) = 0 \, \forall u \in F\,\}$. 
\end{proof}
\begin{lemma}
$F^{\perp}$ is a closed vector subspace of $E_{\omega}$ for all 
vector subspaces $F$ of $\eo$.
\end{lemma} 

\begin{proof}
Let $v_{n}$ be a sequence of $F^{\perp}$ that converges in $\eo$ to 
some $v \in \eo$, i.e. $v = \lim_{n \rightarrow \infty} (v_{n})$.
Then we have that $\forall u \in F, (u,v) = \lim_{n\rightarrow 
infty} (u,v_{n})$. $(u,v_{n}) = 0\, \forall n\in \mathbb{N}$, since
$v_{n} in F^{\perp}$. This implies $(u,v) = 0, \Rightarrow v \in 
F^{\perp}$. Therefore $F^{\perp}$ is closed in $\eo$. 
\end{proof}

\begin{lemma}
Let $\lambda$ be an eigenvalue of $D$ with finite multiplicity with
respect to $D$. Then $R(\lid) = N(\lid)^{\perp}$, therefore 
$R(\lid)$ is closed in $\eo$. 
\end{lemma}
\begin{proof}
We have shown in Theorem 1 (2) that $R(\lid)$ is generated by 
$E_{k}$ such that $\lambda_{k} \ne \lambda$ and in the first part of
Theorem 1 (1) that $N(\lid)$ is generated by the 
$e_{i}$ such that $\lambda_{l} = \lambda$. 
\begin{enumerate}
\item[($1^{\circ}$)] Let $v \in R(\lid)$ with $v = \sum v_{i}e_{i}$
and $v_{i} = 0 $ if $\lambda_{i} = \lambda$. Let $u \in 
N(\lid)$ with $u = \sum u_{j}e_{j}$ with $u_{j} = 0$ if $\lambda_{j}
\ne \lambda$. $(u,v) = \sum u_{k}v_{k}$ if $\lambda_{k}= \lambda
\Rightarrow v_{k}=0$, and if $\lambda_{k} \ne \lambda$, $u_{k} = 0
\Rightarrow (u,v) = 0 \Rightarrow v \in N(\lid)^{\perp}$.
Therefore $R(\lid) \subseteq N(\lid)^{\perp}$.

\item[($2^{\circ}$)]  Let $v \in N(\lid)^{perp}$, $v = \sum_{i \in 
\mathbb{N}} v_{i}e_{i}$. Let $e_{j}$ be such that $\lambda_{j} =
\lambda \Rightarrow e_{j} \in N(\lid)$. Then we have
\begin{displaymath}
\begin{aligned}
(e_{j}, v) &= 0.\, \textrm{ Since } \\
(e_{j},v) &= \sum v_{i}(e_{j}, e_{i}) \, \textrm{ and }\\
(e_{j},e_{i}) &= \left\{\begin{array}{cc} 
0 & \text{ if } i\ne j \\
1 & \text{ if } i = j \\
\end{array}
\right. \\
\Rightarrow (e_{j},v)&= v_{j} \Rightarrow v_{j} = 0
\end{aligned}
\end{displaymath}
Therefore $v$ belongs to the subspace generated by the $e_{k}$ 
such that $\lambda_{k} \ne \lambda$. Therefore we have
$N(\lid)^{\perp} \subseteq R(\lid)$. Therefore $N(\lid)^{\perp}
= R(\lid) \Rightarrow R(\lid)$ is closed in $\eo$. 

Since $\eta(\lid)$ and $\delta(\lid)$ are finite such that
$\eta(\lid) = \delta(\lid) = \#\{e_{i}\,|\,\lambda_{i} = \lambda
\} \Rightarrow X(\lid) = 0$. Since $R(\lid)$ is closed, $(\lid)$
is a Fredholm operator of index 0. 
\end{enumerate}
\end{proof}

\begin{defin}
We define the \emph{proper spectrum} of $D$ to be the set of all
$\lambda_{k} \in \sigma(D)$ with multiplicity 0 with respect to $D$ 
denoted by
\begin{displaymath}
\begin{aligned}
\sigma_{e}^{\prime}(D) &= \left\{ \lambda \in \kf,\, (\lid) 
\text{ is invective but not surjective }\right\} \\
N(\lid) = \{\,0\,\} \\
\end{aligned}
\end{displaymath}
\end{defin}
Note that if $\lambda \in \sigma(D)$ and $\lambda \not \in 
\sigma_{p}(D) = \Lambda$, $\lambda$ is not an eigenvalue $\Rightarrow
N(\lid) = \{\,0\,\}\Rightarrow  \lambda $ is of multiplicity 0 
with respect to $D$. 

\begin{defin}
We define the \emph{improper essential} spectrum of $D$ to be the st of 
all $\lambda_{k} (k \in \mathbb{N})$ with infinite multiplicity with 
respect to $D$, denoted by 
\begin{displaymath}
\sigma_{e}^{\prime \prime}(D) = \left\{ \lambda \in \kf,\,
\eta(\lid) \text{ is infinite }\right\}
\end{displaymath}
where $\eta(\lid) = dim N(\lid)$. 
\end{defin}

\begin{thm}
Let $D: \eo  \longrightarrow \eo$ be a bounded diagonal operator such 
that $D(u) = \sum \lambda_{i}u_{i}e_{i}$ with $u = \sum u_{i}e_{i}$. 
Then $\sigma_{e}(D) = \sigma_{e}^{\prime}(D) \cup \sigma_{e}^{\prime
\prime}(D)$.
\end{thm}

\begin{proof}
\begin{enumerate}
\item{} By the definition of $\sepp(D)$, we conclude that 
$\sepp(D) \subseteq \sigma_{e}(D)$, since $\sigma_{e}(D)$ contains
all the $\lambda$ with infinite multiplicity with respect to $D$ and
all the $\lambda$ with multiplicity $O$ that are not surjective
with respect to $D$. 

\item{} $\lambda \in \sep(D) \Rightarrow$ $\lambda \in \sigma_{e}(D)$
is of multiplicity $0$ with respect  to $D$. This implies that 
$\lambda \not \in  \sigma_{p}(D) \Rightarrow \lambda$ is not an 
eigenvalue. But we know that $\sigma(D) = \sigma_{p}(D) \cup
\sigma_{e}(D) \Rightarrow \lambda \in \sigma_{e}(D) \Rightarrow 
\sigma^{\prime}(D) \subseteq \sigma_{e}(D)$. Therefore 
$\sep(D) \cup \sepp(D) \subseteq \sigma_{e}(D)$. 

\item{} Let $\lambda \in \sigma_{e}(D)$. There are two cases:

Case 1: Let $\lambda \in \sigma_{e}(D)$. Then $\lid$ is not a Fredholm
operator of index $0 \Rightarrow \lambda$ is of infinite multiplicity
$\Leftrightarrow \lambda \in \sepp(D)$.

Case 2: $\lambda \not \in \sigma_{p}(D) = \Lambda \Rightarrow 
N(\lid) = \{\,0\,\}$ and $\lid$ is not surjective. Then $\lambda$ is 
an element of $\sigma(D)$ with multiplicity $0$. Then $\lambda \in 
\sep(D)$. 

\begin{displaymath}
\begin{aligned}
\therefore \sigma(D) &\subseteq \sep(D) \cup \sepp(D)\\
\textrm{ Hence } \sigma_{e}(D) & \subseteq \sep(D) \cup \sepp(D)\\
\end{aligned}
\end{displaymath}
\end{enumerate}
\end{proof}

\begin{thm}
Let $D : \eo \rightarrow \eo$ be a bounded diagonal operator, i.e.
$D(u) = \sum \lambda_{i} u_{i} e_{i}$ with $u = \sum u_{i}e_{i}$ for 
$\lambda_{i} \in \kf$ $\forall i \in \mathbb{N}$. Let $\Lambda =
\{\,\lambda_{k}\,\}_{k\in \mathbb{N}}$ -- $\lambda_{k}$ are the 
eigenvalues. Then $\sep(D) = \partial \Lambda = \overline{\Lambda}
- \Lambda$.  
\end{thm}

\begin{proof}
\begin{enumerate}
\item{}Let $\lambda \in \sep(D)$. Then $\lambda \in \sigma(D)$ and $\lambda$ 
is of multiplicity $0$ with respect to $D$. This implies that
$\lambda \not \in \sigma_{p}(D) = \Lambda$ or $\sigma(D) = 
\overline(\Lambda)$. Thus $\lambda \in \overline{\Lambda} - \Lambda
= \partial \Lambda$. Therefore $\sep(D) = \partial \Lambda$.

\item{}Now let $\lambda \in \plam \Rightarrow \lambda \not \in \Lambda =
\sigma_{p(D)}$. This implies $\lambda$ is of multiplicity $0$ with
respect to $D$ and $\lambda \in \sigma(D) = \overline{\Lambda}
\supseteq \plam \Rightarrow \lambda \in \sep$. Therefore $\plam 
\subseteq \sep(D)$. 
\end{enumerate}
Based on $1$ and $2$, we conclude that
\[
\sep(D) = \plam = \overline{\Lambda} - \Lambda.
\]
\end{proof}

\begin{thm}
Let $D: \eo \rightarrow \eo$ be a bounded diagonal operator. Then 
$\sepp(D) \subseteq \sigma_{p}(D)$.
\end{thm}
\begin{proof}
Every element of $\sepp(D)$ is of infinite multiplicity with respect to 
D $\Leftrightarrow N(\lid) \ne 0 \Rightarrow  \lambda \in \spp(D)$. 
Therefore $\sepp(D) \subseteq \spp(D)$. 
\end{proof}
\begin{cor}
If every eigenvalue of $D$ is of finite multiplicity, then 
$se(D) = \partial \Lambda$. 
\end{cor}
\begin{proof}
Since $\sepp(D) = \phi$ and $ \se(D) = \sep(D) \cup \sepp(D)
\Rightarrow  \sep(D) = \se(D)$
\end{proof}
Note that $\sepp(D) = \phi$ since every eigenvalue of $D$ is of
finite multiplicity and $\sepp(D)$ is the set of all $\lambda$ with 
infinite multiplicity. 
\begin{thm}
Let $D : \eo \rightarrow \eo$ be a linear and diagonal operator of 
finite rank. Then
\begin{enumerate}
\item The spectrum of $D$ is finite, and the spectrum of $D$ equals
$\spp(D)$. 
\item $0$ is an element of the essential spectrum of $D$. 
\end{enumerate}
\end{thm}
\begin{proof}
\begin{enumerate}
\item $D$ is of finite rank, which implies that $\spp(D)$ is
finite. Hence $\se(D)  = \spp(D)$, based on the example (3.2) of DMR

\item Since $\eo$ is infinite dimensional, it implies that $\ker(D)$
is infinite dimensional since $\eo$ is supposed to be of finite rank. This implies that $D$ is not a Fredholm operator (FO). Therefore
$D$ is not a FO of index 0. Thus $0I - D = -D$ is not a 
Fredholm operator of index 0. Therefore $0$ belongs to the essential
spectrum of $D$. 
\end{enumerate}
\end{proof}

\begin{cor}
Let $D$ be a linear diagonal operator of finite rank and let 
$\lambda \in K$ be a nonzero element belonging to $\spp(D)$. Then 
$\lid$ is a Fredholm operator of index $0$. Therefore $\lambda$ does 
not belong to the essential spectrum of $D$. 
\end{cor}
\begin{proof}
The multiplicity of lambda is finite since the spectrum o. Therefore
$\lid$ is a Fredholm operator of index $0$, by Theorem 1. We conclude
that $\lambda$ does not belong to the essential spectrum of $D$. 
\end{proof}

\begin{cor}
If $D$ is a linear diagonal operator of finite rank then the 
essential spectrum of $D$ is equal to $0$.
\end{cor}
\begin{proof}
Since the rank of $D$ is finite, $\ker(D)$ is infinite dimensional. 
This implies that $D$ is not a Fredholm operator, hence $0I-D = -D$
is not a Fredholm operator of index 0. Therefore 0 belongs to the 
essential spectra of $D$. By Theorem 1 and Corollary 1, we conclude 
that the essential spectrum of $D$ is equal to 0. 
\end{proof}

\begin{ex}
Let $\alpha_{n}, n\in \N^{\ast}$ be distinct non zero elements of $\kf$
with $||\alpha_{i}|| < 1$. For each $n\in \N$ let $\beta(n,k) = 
\alpha_{n}$ with $k \in \N$. Now set $\lambda_{i} = \beta(\varphi_{2}(i))$ where $\varphi_{2}: \N \rightarrow \N \times \N$ is a bijection. 
Let $D : \eo \rightarrow \eo$ be defined by  $D(u) = \sum_{i \in \N}
\lambda_{i}u_{i}e_{i}$ for $u = \sum_{i \in \N}u_{i}e_{i} \in \eo$. 

For $i_{0} \in \N$, if $\varphi_{2}(i_{0}) = (n_{0},k_{0})
\in \N \times \N$, we have $\lambda_{i_{0}} = \beta(\varphi_{2}(i_{0}))
= \beta(n_{0}, k_{0}) = \alpha_{n_{0}}$. We also have that
$\forall k\in \N, \beta(n_{0}, k) = \alpha_{n_{0}}$. For each $k \in
\N, \left\{\varphi_{2}^{-1}(n_{0}, k)\right\} = B$ is a subset of
$\N$ but $\lambda_{j_{k}}= \lambda_{i_{0}} = \alpha_{n_{0}}$ where
$j_{k} = \varphi_{2}^{-1}(n_{0},k)$. This implies that 
$N(\lambda_{i_{0}}I - D)$ for all $k \in \N$. This implies 
$D(e_{j_{k}}) = \lambda_{j_{k}}e_{k} \Rightarrow D(e_{j_{k}}) =
\lambda_{i_{0}}e_{k}$. Then $(\lambda_{i_{0}}I - D)(e_{j_{k}}) = 0$.

Since the kernel of a Fredholm operator is finite dimensional, 
$\lambda_{i_{0}}I - D$ is not a Fredholm operator of index 0. This means 
that $\lambda_{i_{0}} \in  \se(D)$ $\forall i_{0}\in \N$. $\therefore
\Lambda \subseteq \se(D)$ since $\partial \Lambda \subseteq 
\se(D)$ by Theorem 2, which implies that
$\se(D) \supseteq \Lambda \cup \partial \Lambda = \overline{\Lambda}
= \sigma(D)$. But $\sigma(D) = \spp(D) \cup \se(D) \Rightarrow
\se(D) \subseteq \sigma(D)$. $\therefore \se(D) = \sigma(D)$. 
\end{ex}

\end{document}